# Polynomial Transformation Method for Non-Gaussian Noise Environment

Jugalkishore K. Banoth, Pradip Sircar*

Department of Electrical Engineering
Indian Institute of Technology Kanpur
Kanpur 208016, India
*Corresponding author; Email address: sircar@iitk.ac.in

**Abstract** - *Signal processing in non-Gaussian noise environment is addressed in this paper. For many real-life situations, the additive noise process present in the system is found to be dominantly non-Gaussian. The problem of detection and estimation of signals corrupted with non-Gaussian noise is difficult to track mathematically. In this paper, we present a novel approach for optimal detection and estimation of signals in non-Gaussian noise. It is demonstrated that preprocessing of data by the orthogonal polynomial approximation together with the minimum error-variance criterion converts an additive non-Gaussian noise process into an approximation-error process which is close to Gaussian. The Monte Carlo simulations are presented to test the Gaussian hypothesis based on the bicoherence of a sequence. The histogram test and the kurtosis test are carried out to verify the Gaussian hypothesis.*

## 1   Introduction

In the signal detection and estimation problems, we often assume that the additive random noise process is Gaussian distributed because this distribution is simple and mathematically tractable, and the assumption makes analytical results possible. However for many real-life situations, the additive noise process is found to be dominantly non-Gaussian. Some examples are the ocean acoustic noise and the urban radio-frequency (RF) noise [1]. The RF receivers designed to perform in white Gaussian noise can not perform satisfactorily when the electromagnetic environment encountered by the receiver system is non-Gaussian in nature [2]. For detection and estimation of radar signals in high clutter environments and similar processing of sonar signals in presence of high reverberation, we need to deal with non-Gaussian noise [1, 2].

There are two existing approaches for solving the problems of detection and estimation of signals in non-Gaussian noise environment. The first approach is to use the robust statistics in lieu of the classical mathematical statistics, and to look for procedures which are consistent or in other words, insensitive to deviations of the noise distribution from the idealized model, i.e., the Gaussian distribution [3]. An optimally robust procedure minimizes the maximum degradation of performance for a preset deviation of the noise distribution. The robust techniques, however, can not provide consistent performance for a noise process with an arbitrary probability density function (PDF).

The second approach to deal with a non-Gaussian noise environment is to use a noise model which is general enough to depict an arbitrary PDF, yet the model retains the desirable simplicity of manipulation as that of a Gaussian PDF. Accordingly, the Gaussian-mixture PDF, the generalized Gaussian PDF, the Middleton class A PDF, and some such PDFs are employed to model non-Gaussian noise [4]. Incidentally, as the noise model is required to be more accurate, the ease of analysis as that of a Gaussian PDF disappears.

In this paper, we present a third approach to deal with a non-Gaussian noise environment, by employing the polynomial transformation method. Preprocessing of data by the orthogonal polynomial approximation (OPA) together with the minimum error-variance criterion (MEC) has an excellent noise-rejection capability [5, 6]. The OPA based transformation was originally proposed to convert non-uniformly sampled data into uniformly sampled data [5]. However, since the transformation provides significant signal enhancement by rejecting the high frequency interference, preprocessing of data may be useful in detection and estimation problems for better accuracy even for uniformly sampled data [6]. Perhaps the most desirable feature of preprocessing the signal samples by the OPA based method is that the statistical distribution of the approximation-error process in the preprocessed data becomes close to Gaussian when the noise process is not necessarily Gaussian distributed [5]. Based on this argument, the maximum likelihood estimator (MLE) can be designed to estimate parameters of a signal corrupted with non-Gaussian noise [7].

In the present work, we take a closer look of the preprocessing of data by the OPA based method, and we test the hypothesis that the approximation-error process in the preprocessed data is Gaussian distributed even when the noise process corrupting the sampled data is non-Gaussian. Several types of tests are applied for testing the hypothesis. We plot the histogram of a given sequence and look for the proverbial bell shape as a simple test for its Gaussian distribution [8]. We compute the kurtosis [9] and apply the Hinich test [10, 11] for validation of the Gaussian hypothesis. We consider the following noise processes for the Monte Carlo simulation: (i) Gaussian, (ii) Laplacian, (iii) Uniform, and (iv) Gamma distributed.

## 2 Orthogonal Polynomial Transformations

The real-valued discrete-time signal $g[n]$ is to be detected/estimated utilizing the sampled sequence $x[n] = g[n] + w[n]$, where $w[n]$ is the noise sequence which may not be Gaussian distributed. The sampled data $\{x[n]\}$ are preprocessed by the orthogonal polynomial transformation to obtain the transformed data $\{y[n]\}$ as follows [6],

$$\mathbf{y} = \mathbf{P}\mathbf{Q}^{-1}\mathbf{P}^T\mathbf{x} \quad (1)$$

where

$\mathbf{x} = [x[0], x[1], \cdots, x[N-1]]^T$,

$\mathbf{y} = [y[0], y[1], \cdots, y[N-1]]^T$,

$(\mathbf{P})_{ij} = p_j[i]; i = 0, 1, \cdots, N-1; j = 0, 1, \cdots, J-1;$

$$\mathbf{Q} = \mathbf{P}^T\mathbf{P} = \mathrm{Diag}\left[\sum_{m=0}^{N-1}p_0^2[m], \sum_{m=0}^{N-1}p_1^2[m], \cdots, \sum_{m=0}^{N-1}p_{J-1}^2[m]\right]$$

The orthogonal polynomials $p_j[n]$ are computed by the recurrence relation given in [5, 12], and the order of approximation $J$ is chosen such that the error-variance is minimum.

The transformed sequence $y[n]$ is given by $y[n] = g[n] + e[n]$, where $e[n]$ is the approximation error. By utilizing the relation between the error sequence $e[n]$ and the noise sequence $w[n]$,

$$e[n] = \sum_{m=0}^{N-1}\xi_{nm}w[m] \quad (2)$$

where

$$\xi_{nm} = \sum_{j=0}^{J-1}\left\{\frac{p_j[n]p_j[m]}{\sum_{i=0}^{N-1}p_j^2[i]}\right\},$$

we can compute the autocorrelation functions (ACFs) of the error process [6], provided the ACFs of the noise process are known. Furthermore, by invoking the central limit theorem when the random variables $w[n]$ are independent with zero mean and identical variance, and the coefficients $\xi_{nm}$ are bounded [5, 13], we can argue that the error process will be close to Gaussian even when the distribution of the noise process is non-Gaussian.

## 3 Gaussian Hypothesis Testing

The third order cumulant of the noise/ error process $u[n]$ is given by

$$C_{3u}[k_1, k_2] = E\{u[n]u[n+k_1]u[n+k_2]\} \quad (3)$$

where $E$ is the expectation operator, and the third order spectrum, commonly known as the bispectrum, is defined as the two-dimensional Fourier transform of the third order cumulant,

$$S_{3u}(\omega_1, \omega_2) = \sum_{k_1=-\infty}^{\infty}\sum_{k_2=-\infty}^{\infty}C_{3u}[k_1, k_2]\exp\{-j(\omega_1 k_1 + \omega_2 k_2)\} \quad (4)$$

The squared bicoherence $|B_{3u}(\omega_1, \omega_2)|^2$ is determined as follows

$$|B_u(\omega_1, \omega_2)|^2 = \frac{|S_{3u}(\omega_1, \omega_2)|^2}{S_{2u}(\omega_1)S_{2u}(\omega_2)S_{2u}(\omega_1+\omega_2)} \quad (5)$$

where $S_{2u}(\omega)$ is the power spectrum.

The Hinich test is based upon the squared bicoherence at a bifrequency $(\omega_1, \omega_2)$ being zero for Gaussianity of the underlying sequence. The $|B|^2$ value is averaged over the principal domain [10, 11], and the resulting statistics is central $\chi^2$ distributed under the null hypothesis: $S_{3u}(\omega_1, \omega_2) \equiv 0$. Hence, it is easy to devise a statistical test to determine whether the observed squared bicoherence is consistent with a central $\chi^2$ distribution by computing a probability of false alarm (PFA) value. If the null hypothesis of the bispectrum being zero is not rejected, we then compute the average kurtosis $K_u$ given by [9]

$$K_u = \frac{E\{u^4\}}{\left[E\{u^2\}\right]^2} - 3 \qquad (6)$$

where the value is averaged over each element of the sequence $\{u[n]\}$. The kurtosis test is based on the null hypothesis: $K_u \equiv 0$ for a Gaussian distribution of the underlying process.

## 4 Simulation Results

The real part of the complex-exponential transient discrete-time signal

$$g[n] = \sum_{i=1}^{3} b_i \exp(j\varphi_i)\exp(s_i nT) \qquad (7)$$

where
$b_1 = 1.0, s_1 = -0.2 + j2.0, \varphi_1 = 0;$
$b_2 = 0.5, s_2 = -0.1 + j4.0, \varphi_2 = \pi/4;$
$b_3 = 0.5, s_3 = -0.3 + j1.0, \varphi_3 = \pi/6;$

corrupted with non-Gaussian noise setting the signal-to-noise ratio (SNR) at 10 dB, is sampled at 60 uniformly spaced points with time interval $T = 0.15$.

We present the three cases of Laplacian, Uniform, and Gamma distributed noise environments in this work, beside the Gaussian noise case. In each case, after applying the polynomial transformation we obtain the transformed data, and then, subtracting the transient signal from the transformed data, the error process is separated. The input noise and the output error processes are tested for Gaussianity. Figs. 1−4 show the bispectrum and the histogram plots.

For the Gaussian noise case, we calculate the bicoherence of the output error and check whether the squared bicoherence is consistent with a central $\chi^2$ distribution by computing the PFA value. The PFA is computed to be 0.9479, which is high, and we cannot reject the null hypothesis. The average kurtosis value for the output error is computed to be −0.1526, whereas the kurtosis value for the input noise is computed to be −0.0857 (theoretical value zero). For the Laplacian case, the PFA for the input noise is 0.396, and the PFA for the output error is 0.9975. Since the PFA of the output error is high, we cannot reject the null hypothesis. The kurtosis values are 2.9359 for the input noise and 0.0148 for the output error. For the Uniform noise case, we compute the PFA for the input noise to be 0.6973 and the PFA for the output error to be 0.9649. The average kurtosis values are computed to be −1.2408 for the input noise and −0.1346 for the output error. For the Gamma distributed noise environment, the PFA for the input noise is 0.7379, and the PFA for the output error is 0.9847. The kurtosis values are 0.6962 for the input noise and 0.1025 for the output error. In all cases, we find that the average kurtosis value of the output error process is near zero, confirming that the error process is close to Gaussian.

## 5 Concluding Remarks

In this paper, we present a new technique for optimal detection and estimation of signals corrupted with non-Gaussian noise. We preprocess the sampled data by the polynomial transformation method which converts the noise process into an approximation-error process which is Gaussian distributed.

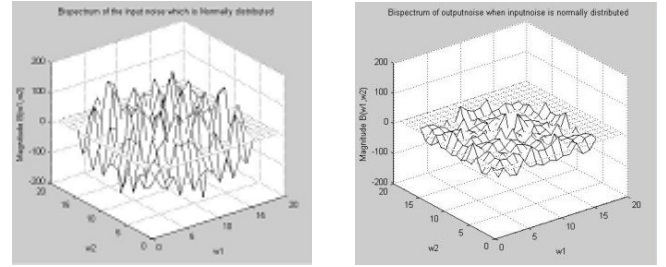

Figure 1(a): The Bispectrum of the input noise (Gaussian) and the output error process

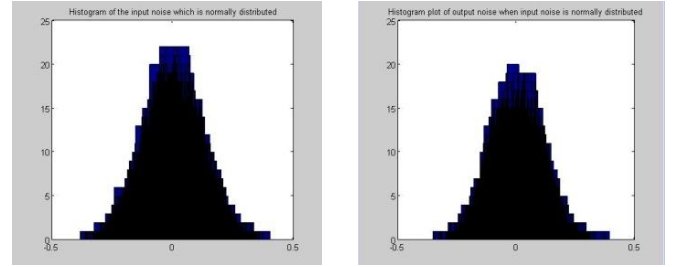

Figure 1(b): The Histogram of the input noise (Gaussian) and the output error process

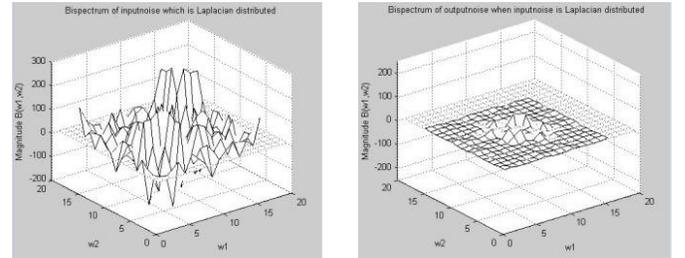

Figure 2(a): The Bispectrum of the input noise (Laplacian) and the output error process

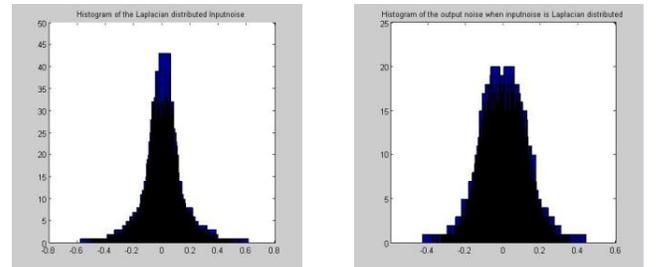

Figure 2(b): The Histogram of the input noise (Laplacian) and the output error process

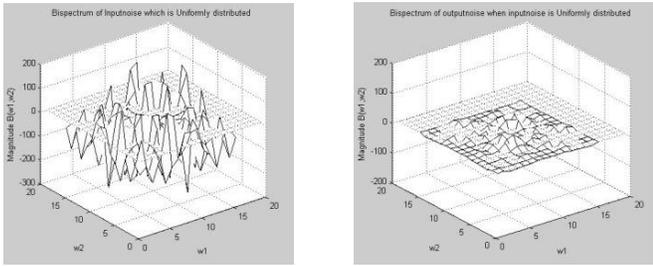

Figure 3(a): The Bispectrum of the input noise (Uniform) and the output error

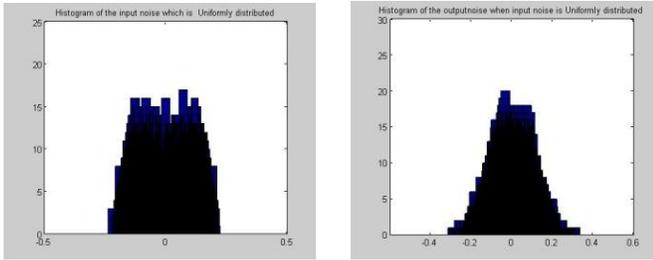

Figure 3(b): The Histogram of the input noise (Uniform) and the output error

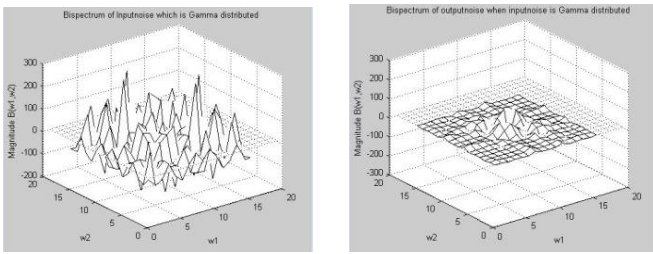

Figure 4(a): The Bispectrum of the input noise (Gamma distributed) and the output error

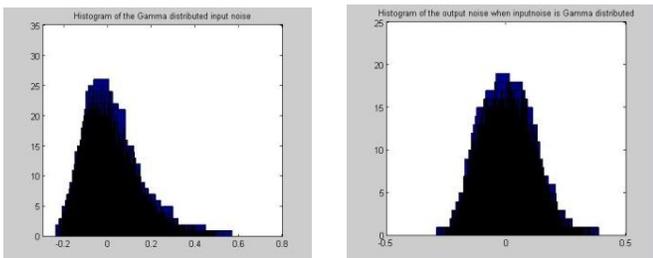

Figure 4(b): The Histogram of the input noise (Gamma distributed) and the output error